\documentclass[a4paper,fleqn]{article}
\usepackage{amssymb,amsthm,amsfonts}
\usepackage{amsmath}
\newtheorem{theorem}{Theorem}

\newtheorem{corollary}{Corollary}

\newtheorem{remark}{Remark}

\numberwithin{equation}{section}
\numberwithin{lemma}{section}
\numberwithin{theorem}{section}
\numberwithin{corollary}{section}

\allowdisplaybreaks
\begin{document}
\setcounter{page}{1}

\title{Some New Results for Generalized  Incomplete Exponential  Matrix Functions}

\author{Ashish Verma
\\ 
Department of Mathematics\\ Prof. Rajendra Singh (Rajju Bhaiya)\\ Institute of Physical Sciences for Study and Research \\  V. B. S. Purvanchal University, Jaunpur  (U.P.)- 222003, India\\
vashish.lu@gmail.com (Corresponding author)
\\[10pt]
Komal Singh Yadav
\\ 
Department of Mathematics\\ Prof. Rajendra Singh (Rajju Bhaiya)\\ Institute of Physical Sciences for Study and Research \\  V. B. S. Purvanchal University, Jaunpur  (U.P.)- 222003, India\\
ksyadav230128@gmail.com\\[10pt]
}

\maketitle
\begin{abstract}

The primary goal of this paper is to introduce and investigate generalized incomplete exponential functions with matrix parameters. Integral representation, differential formula, addition formula, multiplication formula, and recurrence relation obtained here are believed to be new in the theory of special matrix functions. We also establish the connection between these matrix functions and other matrix functions, such as the incomplete gamma matrix function, the Bessel and modified Bessel matrix functions.\\[12pt]

\end{abstract}

\section{Introduction}
Tricomi  studied the theory of incomplete gamma functions \cite{FG}. These are very important special functions that are used in a variety of problems in mathematics, astrophysics, applied statistics, and engineering. These functions are also useful in the study of   Fourier and Laplace transforms, as well as probability theory. The incomplete exponential functions  and the  incomplete hypergeometric functions have been introduced  by Pal {\em et al}. \cite{APa, AQ}. They also presented applications in communication theory, probability theory, and groundwater pumping modelling. The incomplete Pochhammer symbols  and the  incomplete hypergeometric functions have been introduced  by Srivastava {\em et al}. \cite{HM}.

The matrix theory is used in orthogonal polynomial and special functions, and it is widely used in mathematics in general. The special matrix functions are mentioned in the literatures of statistics \cite{A} and Lie theory \cite{AT}. It is also related to the matrix version of the Laguerre, Hermite, and Legendre differential equations, and the corresponding polynomial families are listed in \cite{LR,LRE,LR1}. J\'odar and  Cort\'es  has  introduced  matrix analogue of the Gauss hypergeometric function in \cite{LC}. Abdalla \cite{Ab} introduced the incomplete hypergeometric matrix functions  and  discussed some fundamental properties of these functions.  The Wright hypergeometric matrix functions and the incomplete Wright Gauss hypergeometric matrix functions have been established in \cite{AY}  and  obtained some properties of these functions.   In \cite{ A3, AP},  introduced the incomplete second Appell hypergeometric matrix functions and the incomplete Srivastava's triple hypergeometric matrix functions and studied some basic properties of these functions, including matrix differential equation, integral formula, recursion formula, recurrence relation, and differentiation formula.

The paper is organized in the following manner.  In Section~2, we list the fundamental definitions and results of special matrix functions that are required in the sequel.  In Section~3, we define incomplete exponential matrix functions $e(x;t)$ and $E(x;t)$. Some properties;  integral representation,  differentitation formula and connections with bessel matrix functions are also derived. In Section~4, we introduce  incomplete exponential  matrix functions $_{p}e_{q}(x;t)$ and $ _{p}E_{q}(x;t)$ and investigate several properties of each of these functions. Finally, in Section 5, we define generalized incomplete exponential  matrix functions $_{p}e_{q}(x, A, B;v)$ and $ _{p}E_{q}(x, A, B;v)$ and derive certain properties of each of these  functions. 

\section{ Preliminaries}
Throughout this paper,  let $\mathbb{C}^{r\times r}$ be the vector space of $r$-square matrices with complex entries. For any matrix $E\in \mathbb{C}^{r\times r}$, its spectrum $\sigma(E)$ is the set of eigenvalues of $E$. 
 A square matrix $E$ in $\mathbb{C}^{r\times r}$  is said to be positive stable if $\Re(\lambda)>0$ for all $\lambda\in\sigma(E)$. 

Let $E$ be positive stable matrix in $\mathbb{C}^{r\times r}$. The gamma matrix function $\Gamma(E)$ are defined as follows \cite{LC1}:
\begin{align}\Gamma(E)=\int_{0}^{\infty} e^{-t} t^{E-I} dt; \hskip1cm t^{E-I}= \exp((E-I)\ln t). 
\end{align}

 The reciprocal gamma function \cite{LC1} is defined as 
$\Gamma^{-1}(E)= (E)_n \ \Gamma^{-1}(E+nI)$, where  $E+nI$ is invertible for all integers $n\geq 0$, $I$ being the $r$-square identity matrix and $(E)_n$ is the shifted factorial matrix function for $E\in\mathbb{C}^{r\times r}$ defined in \cite{LC}.

If $E\in\mathbb{C}^{r\times r}$ is a positive stable matrix and $n\geq1$, then by \cite{LC1} we have $\Gamma(E) = \lim_{n \to\infty}(n-1)! (E)^{-1}_{n} n^E$.  

By  application of the matrix functional calculus, the Pochhammer symbol  for  $E\in \mathbb{C}^{r\times r}$ is given by \cite{LC1}
\begin{equation}
(E)_n = \begin{cases}
I, & \text{if $n = 0$},\\
E(E+I) \dots (E+(n-1)I), & \text{if $n\geq 1$}.
\end{cases}\label{c1eq.09}
\end{equation}
This gives
\begin{equation}
(E)_n = \Gamma^{-1}(E) \ \Gamma (E+nI), \qquad n\geq 1.\label{c1eq.010}
\end{equation} 

If $E$ and $F$ are positive stable matrices in $\mathbb{C}^{r \times r}$, then, for $EF = FE$, the beta matrix function is defined as \cite{LC1}
\begin{align}
\mathfrak{B}(E,F) =\Gamma(E)\Gamma(F)\Gamma^{-1}(E+F) &= \int_{0}^{1}t^{E-I}(1-t)^{F-I}dt \label{1ca1.4}\\
&=\int_{0}^{\infty}u^{E-I}(1+u)^{-(E+F)}du.\label{1ca1.5}
\end{align}
Clearly, the generalized pochhammer matrix symbol $(A)_{kn}$  can be represented in the following form
\begin{align}
(E)_{kn}= k^{kn}\,\left (\frac{E}{k}\right)_{n}  \left(\frac{E+I}{k}\right)_{n}\cdots  \left(\frac{E+(k-1)I}{k}\right)_{n}.\label{g1}
\end{align}
Suppose $E$ be positive stable matrix in $\mathbb{C}^{r\times r}$ and let $x$ be a positive real number.
Then the  incomplete gamma matrix  functions $\gamma(E,x)$ and $\Gamma(E,x)$ are defined by \cite{Ab}
\begin{align}\gamma(E,x)=\int_{0}^{x} e^{-t} t^{E-I}dt \label{1eq4}
\end{align}
and
\begin{align}
 \Gamma(E,x)= \int_{x}^{\infty} e^{-t} t^{E-I}dt\,,\label{1eq5}
\end{align}
respectively and satisfy the following decomposition formula:
\begin{align}\gamma(E,x)+\Gamma(E,x)=\Gamma(E).\label{1eq6}
\end{align}
Let  $E$  be matrix in $\mathbb{C}^{r\times r}$ and let x be a positive real number. Then the incomplete Pochhammer matrix symbols $(E;x)_{n}$ and $[E; x]_{n}$ are defined as follows  \cite{Ab}
\begin{align}
(E; x)_{n}= \gamma(E+nI, x) \,\Gamma^{-1}(E)\label{1eq7}
\end{align}
and 
\begin{align}
[E; x]_{n}= \Gamma(E+nI, x)\, \Gamma^{-1}(E).\label{1eq8}
\end{align}
In idea of (\ref{1eq6}), these incomplete Pochhammer matrix symbols $(E; x)_{n}$ and $[E; x]_{n}$ complete the following decomposition relation
\begin{align}(E; x)_{n}+[E; x]_{n}= (E)_{n}.\label{1eq9}
\end{align}
where $(E)_{n}$ is the Pochhammer matrix  symbol defined in \cite{LC}.

Let $E$, $F$ and $G$ are matrices in  $\mathbb{C}^{r\times r}$ such that $G+nI$ is invertible for all integers $n \geq0.$ The incomplete Gauss hypergeometric matrix  functions  are defined by  \cite{Ab}
\begin{align}{ _2\gamma_{1}}\Big[(E; x), F; G; z\Big]= \sum_{n=0}^{\infty}(E;x)_{n} (F)_{n}(G)_{n}^{-1}\frac{z^{n}}{n!}\label{11eq10}
\end{align}
and
\begin{align}{_2\Gamma_{1}}\Big[[E; x], F; G; z\Big]= \sum_{n=0}^{\infty}[E;x]_{n} (F)_{n}(G)_{n}^{-1}\frac{z^{n}}{n!}.\label{1eq10}
\end{align}

The matrix function $_{p}R_{q}(A, B; z)$ is  \cite{RR} defined by:
\begin{align}
_{p}R_{q}(A, B; z)&= {_{p}R_{q}}\left[\begin{array}{c}E_1, E_2, \cdots, E_{p}\\ F_1, F_2, \cdots, F_{q}\end{array}\vert A, B; v \right]\notag\\
&= \sum_{m=0}^{\infty}\Gamma^{-1}(mA+B) (E_1)_m\cdots (E_p)_{m}(F_1)_m\cdots (F_q)_{m}\frac{v^{m}}{m!},\label{9}
\end{align}
where $A$, $B$, $E_i$ and $F_{j}$, $1\leq i\leq  p$, $1\leq j\leq  q$ be matrices in $\mathbb{C}^{r\times r}$ such that $F_{j}+kI$ are invertible for all integers $k\geq 0$.

The  Bessel matrix function is \cite{J1, J2, J3} defined by: 
\begin{align}J_{A}(z)=\sum_{m\geq 0}^{\infty}\frac{(-1)^{m}\,\,\Gamma^{-1}(A+(m+1)I)}{m!}\Big(\frac{z}{2}\Big)^{A+2mI},\label{j1}
\end{align}
where $A+nI$ is invertible for all integers $n\geq 0$. Therefore, the modified Bessel matrix functions are introduced in \cite{J3}  in the form 
\begin{align}&I_{A}= e^\frac{-Ai\pi}{2} J_{A}(z  e^\frac{i\pi}{2}); \,\,\, -\pi<arg(z)<\frac{\pi}{2},\notag\\
&I_{A}= e^\frac{Ai\pi}{2} J_{A}(z  e^\frac{-i\pi}{2}); \,\,\, -\frac{\pi}{2}<arg(z)<\pi.\label{j2}
\end{align}

\section{Incomplete Exponential Matrix Functions $e(x;t)$ and $E(x;t)$: Definition and some new properties}
In this section, we  define the incomplete exponential matrix functions as 
\begin{align}
e\left((x, t); A\right)= \sum_{m=0}^{\infty} {\Gamma^{-1}(A+mI)}{\gamma(A+mI, x)}\frac{t^{m}}{m!},\label{e21}
\end{align}
\begin{align}
E\left((x, t); A\right)= \sum_{m=0}^{\infty}{\Gamma^{-1}(A+mI)}{\Gamma(A+mI, x)}\frac{t^{m}}{m!},\label{e22}
\end{align}
where $A$ is a positive stable matrix in $\mathbb{C}^{r\times r}$ such that $A+kI$  is  invertible for all integers $k\geq0$. 
So that
\begin{align}
e\left((x, t); A\right)+ E\left((x, t); A\right)= e^{t}.\label{e23}
\end{align}
\begin{theorem}
 Let $A$ be positive stable matrix in  $\mathbb{C}^{r\times r}$ such that $A+kI$ is  invertible for all integers $k\geq0$. Then the following integral representation for the incomplete exponential matrix functions  holds true:
\begin{align}
e\left((x, t); A\right)&= \Gamma^{-1}(A)\int_{0}^{x} v^{A-I} e^{-v}\Big(\sum_{m=0}^{\infty}(A)^{-1}_{m}\frac{(vt)^{m}}{m!}\Big) dv\notag\\
&= \Gamma^{-1}(A)\int_{0}^{x} v^{A-I} e^{-v}\, _{0}F_{1}(-; A; vt)dv,\label{e24}
\end{align}
\begin{align}
E\left((x, t); A\right)&= \Gamma^{-1}(A)\int_{x}^{\infty} v^{A-I} e^{-v}\Big(\sum_{m=0}^{\infty}(A)^{-1}_{m}\frac{(vt)^{m}}{m!}\Big) dv\notag\\
&= \Gamma^{-1}(A)\int_{x}^{\infty} v^{A-I} e^{-v}\, _{0}F_{1}(-; A; vt)dv.\label{e25}
\end{align}
\end{theorem}
\begin{proof} From definition of  the incomplete exponential matrix functions, we have
 \begin{align}
e\left((x, t); A\right)= \sum_{m=0}^{\infty} {\Gamma^{-1}(A+mI)}{\gamma(A+mI, x)}\frac{t^{m}}{m!}.
\end{align}
Applying the definition of incomplete gamma matrix  functions (\ref{1eq4}),  we get (\ref{e24}). In a similar manner we can prove (\ref{e25}). This finishes the proof of this theorem.
\end{proof}
\begin{corollary}(Connections with Bessel matrix functions )The following integral representation for the incomplete exponential matrix functions  holds true:
\begin{align}
&e\left((x, t); A+I\right)= t^{-\frac{A}{2}}\int_{0}^{x} v^{\frac{A}{2}} e^{-v} I_{A}(2\sqrt {vt}) dv,\\
&E\left((x, t); A+I\right)= t^{-\frac{A}{2}}\int_{x}^{\infty} v^{\frac{A}{2}} e^{-v} I_{A}(2\sqrt {vt}) dv,\\
&e\left((x, -t); A+I\right)= t^{-\frac{A}{2}}\int_{0}^{x} v^{\frac{A}{2}} e^{-v} J_{A}(2\sqrt {vt}) dv,\\
&E\left((x, -t); A+I\right)= t^{-\frac{A}{2}}\int_{x}^{\infty} v^{\frac{A}{2}} e^{-v} J_{A}(2\sqrt {vt}) dv,
\end{align}
where $I_{A}(v)$ and $J_{A}(v)$ are Bessel matrix functions defined in \cite{J3}.
\end{corollary}
\begin{theorem}
 Let $A$ be positive stable matrix in  $\mathbb{C}^{r\times r}$ such that $A+kI$ is  invertible for all integers $k\geq0$.  Then the following derivative formula  for the incomplete exponential matrix functions holds true:
\begin{align}
&\frac{\partial}{\partial t^{n}}e\left((x, t); A\right)= e\left((x, t); A+nI\right),\label{e26}\\
&\frac{\partial}{\partial t^{n}}E\left((x, t); A\right)= E\left((x, t); A+nI\right).\label{e27}\\
&\frac{\partial}{\partial x}e\left((x, t); A\right)= {x^{A-I} e^{-x}}\,{\Gamma^{-1}(A)}\, _{0}F_{1}(-; A; tv),\label{e28}\\
&\frac{\partial}{\partial x}E\left((x, t); A\right)= -{x^{A-I} e^{-x}}\,{\Gamma^{-1}(A)}\, _{0}F_{1}(-; A; tv).\label{e29}
\end{align}
\end{theorem}
\begin{proof}Differentiating   (\ref{e21}) with respect to $t$, the resultant equation comes to:
\begin{align}&\frac{\partial}{\partial t}\,e\left((x, t); A\right)
=\sum_{m=1}^{\infty}{\gamma(A+mI, x)}{\Gamma^{-1}(A+mI)}\frac{t^{m-1}}{(m-1)!}.\label{d4}
\end{align}
 Changing $m$ to $m+1$ and $A$ to $A+I$ in (\ref{d4}) , we obtain
\begin{align}&\frac{\partial}{\partial t}\,e\left((x, t); A\right)= e\left((x, t); A+I\right),
\end{align}
which is (\ref{e26}) for $n=1$. Generalization can be achieved by using the principle of mathematical induction on $n$.

This completes the proof of (\ref{e26}). Successively  (\ref{e27})-(\ref{e29})  can be proved in an analogous manner.
\end{proof}

\section{ Incomplete Exponential  Matrix Functions $_{p}e_{q}(x;t)$ and $ _{p}E_{q}(x;t)$:  Definition and some new relations}
Let $A$, $E_i$ and  $F_j$ , $2\leq i \leq p$, $2\leq j \leq q$, be matrices in  $\mathbb{C}^{r\times r}$ such that  $A+kI$, $F_j+kI$, $2\leq j \leq q$   are invertible for all integers $k\geq0$.  Then we define incomplete generalized exponential  matrix functions as
\begin{align}
_{p}e_{q}\left[(x;t)\vert\begin{array}{c}(A, E_2, \cdots, E_{p})\\ (A, F_2, \cdots, F_{q})\end{array}\right]:=&\sum_{m=0}^{\infty}\, \Gamma^{-1}(A+mI)\gamma(A+mI, x) (E_{2})_{m}\cdots(E_{p})_{m}\notag\\&\times (F_{2})^{-1}_{m}\cdots(F_{p})^{-1}_{m}\frac{t^{m}}{m!},\label{e30}\\
_{p}E_{q}\left[(x;t)\vert\begin{array}{c}(A, E_2, \cdots, E_{p})\\ (A, F_2, \cdots, F_{q})\end{array}\right]:=&\sum_{m=0}^{\infty}\,\Gamma^{-1}(A+mI) \Gamma(A+mI, x) (E_{2})_{m}\cdots(E_{p})_{m}\notag\\&\times (F_{2})^{-1}_{m}\cdots(F_{p})^{-1}_{m}\frac{t^{m}}{m!}.\label{e31}
\end{align}
In view of (\ref{1eq6}), we have the following decomposition formula
\begin{align}
&_{p}e_{q}\left[(x;t)\vert\begin{array}{c}(A, E_2, \cdots, E_{p})\\ (A, F_2, \cdots, F_{q})\end{array}\right]+\, _{p}E_{q}\left[(x;t)\vert\begin{array}{c}(A, E_2, \cdots, E_{p})\\ (A, F_2, \cdots, F_{q})\end{array}\right]\notag\\&=\, _{p-1}F_{q-1}\left[\Big(\begin{array}{c}E_2, \cdots, E_{p}\\  F_2, \cdots, F_{q}\end{array}\Big); t\right]\label{e32}
\end{align}
where $_{p-1}F_{q-1}$ is the generalized hypergeometric matrix function \cite{RD3}.
\begin{theorem}
Let $A$, $E_i$ and  $F_j$ , $2\leq i \leq p$, $2\leq j \leq q$, be matrices in  $\mathbb{C}^{r\times r}$ such that  $A+kI$, $F_j+kI$, $2\leq j \leq q$   are invertible for all integers $k\geq0$. Then the following integral representations for  generalized incomplete exponential  matrix functions holds true:
\begin{align}
&_{p}e_{q}\left[(x;t)\vert\begin{array}{c}(A, E_2, \cdots, E_{p})\\ (A, F_2, \cdots, F_{q})\end{array}\right]\notag\\&= \Gamma^{-1}(A)\int_{0}^{x} v^{A-I} e^{-v}\, _{p-1}F_{q}\left[\Big(\begin{array}{c}-, E_2, \cdots, E_{p}\\ A,  F_2, \cdots, F_{q}\end{array}\Big); vt\right]dv,\label{e33}\\
&_{p}E_{q}\left[(x;t)\vert\begin{array}{c}(A, E_2, \cdots, E_{p})\\ (A, F_2, \cdots, F_{q})\end{array}\right]\notag\\&= \Gamma^{-1}(A)\int_{x}^{\infty} v^{A-I} e^{-v}\, _{p-1}F_{q}\left[\Big(\begin{array}{c}-, E_2, \cdots, E_{p}\\ A,  F_2, \cdots, F_{q}\end{array}\Big); vt\right]dv,\label{e34}\\
& _{p-1}F_{q-1}\left[\Big(\begin{array}{c}E_2, \cdots, E_{p}\\  F_2, \cdots, F_{q}\end{array}\Big); t\right]\notag\\&= \Gamma^{-1}(A)\int_{0}^{\infty} v^{A-I} e^{-v}\, _{p-1}F_{q}\left[\Big(\begin{array}{c}-, E_2, \cdots, E_{p}\\ A,  F_2, \cdots, F_{q}\end{array}\Big); vt\right]dv.\label{e34}
\end{align}
\end{theorem}
\begin{corollary}
The following integral representation for the incomplete   exponential  matrix functions  holds true:
\begin{align}
&_{2}e_{1}\left[(x;t)\vert\begin{array}{c}(C, A)\\ (C)\end{array}\right]= \Gamma^{-1}(C)\int_{0}^{x} v^{C-I} e^{-v}\, _{1}F_{1}\left(A; C; vt\right)dv,\label{e35}\\
&_{2}E_{1}\left[(x;t)\vert\begin{array}{c}(C, A)\\ (C)\end{array}\right]= \Gamma^{-1}(C)\int_{x}^{\infty} v^{C-I} e^{-v}\, _{1}F_{1}\left(A; C; vt\right)dv,\label{e36}\\
& _{2}e_{1}\left[(x;t)\vert\begin{array}{c}(C, A)\\ (C)\end{array}\right]+\, _{2}E_{1}\left[(x;t)\vert\begin{array}{c}(C, A)\\ (C)\end{array}\right]= (1-t)^{-A},\label{e37}
\end{align}
where $A$, $C$  are matrices in  $\mathbb{C}^{r\times r}$ such that $C+kI$ is  invertible for all integers $k\geq0$.
\end{corollary}
\begin{corollary}
The following integral representation for the incomplete  exponential  matrix functions  holds true:
\begin{align}
&_{2}e_{1}\left[(x;-t)\vert\begin{array}{c}(C, C)\\ (C)\end{array}\right]= \Gamma^{-1}(C)\int_{0}^{x} v^{C-I} e^{-(t+1)v}dv,\label{e38}\\
&_{2}E_{1}\left[(x;-t)\vert\begin{array}{c}(C, C)\\ (C)\end{array}\right]= \Gamma^{-1}(C)\int_{x}^{\infty} v^{C-I} e^{-(t+1)v}dv,\label{e39}\end{align}
\end{corollary}where $C$ is a matrices in  $\mathbb{C}^{r\times r}$ such that $C+kI$ is  invertible for all integers $k\geq0$.
\begin{corollary}The following integral representation for the incomplete   exponential  matrix functions  holds true:
\begin{align}
&_{3}e_{1}\left[(x;t)\vert\begin{array}{c}(C, A, B)\\ (C)\end{array}\right]= \Gamma^{-1}(C)\int_{0}^{x} v^{C-I} e^{-v}\, _{2}F_{1}\left(A, B; C; vt\right)dv,\label{e40}\\
&_{3}E_{1}\left[(x;t)\vert\begin{array}{c}(C, A, B)\\ (C)\end{array}\right]= \Gamma^{-1}(C)\int_{x}^{\infty} v^{C-I} e^{-v}\, _{2}F_{1}\left(A, B; C; vt\right)dv,\label{e41}
\end{align}
where $A$, $B$, $C$  are matrices in  $\mathbb{C}^{r\times r}$ such that $C+kI$ is  invertible for all integers $k\geq0$.
\end{corollary}
\section{Generalized Incomplete Exponential  Matrix Functions $_{p}e_{q}(x, A, B;v)$ and $ _{p}E_{q}(x, A, B;v)$:  Definition and some new formulas}
In this section, we introduce the matrix analogs of generalized incomplete exponential functions \cite{APa, AQ}. We denote these matrix functions by $_{p}e_{q}(x, A, B;v)$ and $ _{p}E_{q}(x, A, B;v)$.\\
Let $A$, $B$, $E_i$ and $F_{j}$, $1\leq i\leq  p$, $1\leq j\leq  q$ be matrices in $\mathbb{C}^{r\times r}$ such that $F_{j}+kI$ are invertible for all integers $k\geq 0$. Then, we define 
\begin{align}
_{p}e_{q}(x, A, B;v)\notag&= \, _{p}e_{q}\left[(x, A, B; v)\vert\begin{array}{c}E_1, E_2, \cdots, E_{p}\\ F_1, F_2, \cdots, F_{q}\end{array}\right]\notag\\
&= \sum_{m=0}^{\infty}\Gamma^{-1}(mA+B)\gamma(mA+B, x) (E_{1})_{m}\dots (E_p)_{m} (F_{1})^{-1}_{m}\dots (F_q)^{-1}_{m}\frac{v^{m}}{m!}, \label{e42}
\end{align}
\begin{align}
_{p}E_{q}(x, A, B;v)\notag&= \, _{p}E_{q}\left[(x, A, B; v)\vert\begin{array}{c}E_1, E_2, \cdots, E_{p}\\ F_1, F_2, \cdots, F_{q}\end{array}\right]\notag\\
&= \sum_{m=0}^{\infty}\Gamma^{-1}(mA+B)\Gamma(mA+B, x) (E_{1})_{m}\dots (E_p)_{m} (F_{1})^{-1}_{m}\dots (F_q)^{-1}_{m}\frac{v^{m}}{m!}. \label{e43}
\end{align}
From (\ref{e42}) and (\ref{e43}), we can obtain the following decomposition formula:
\begin{align}
&_{p}e_{q}\left[(x, A, B; v)\vert\begin{array}{c}E_1, E_2, \cdots, E_{p}\\ F_1, F_2, \cdots, F_{q}\end{array}\right]+\, _{p}E_{q}\left[(x, A, B; v)\vert\begin{array}{c}E_1, E_2, \cdots, E_{p}\\ F_1, F_2, \cdots, F_{q}\end{array}\right]\notag\\
&=\, _{p}F_{q}\left[ \begin{array}{c}E_1, E_2, \cdots, E_{p}\\ F_1, F_2, \cdots, F_{q}\end{array}\vert  v\right], \label{e44}
\end{align}
where $_{p}F_{q}$ is the generalized hypergeometric matrix function \cite{RD3}.
\begin{remark}If we take $p=q=0$ and $A=I$,  then (\ref{e42}) and (\ref{e43}) reduces to incopmlete exponential matrix functions  (\ref{e21}) and (\ref{e22}):
\begin{align}
_{0}e_{0}(x, I, B;v)&= \, _{0}e_{0}\left[(x, I, B; v)\vert\begin{array}{c}-\\ -\end{array}\right]\notag\\
&= \sum_{m=0}^{\infty}\Gamma^{-1}(mI+B)\, \gamma(mI+B, x) \frac{v^{m}}{m!}, \label{e45}
\end{align}
and 
\begin{align}
_{0}E_{0}(x, I, B;v)&= \, _{0}E_{0}\left[(x, I, B; v)\vert\begin{array}{c}-\\ -\end{array}\right]\notag\\
&= \sum_{m=0}^{\infty}\Gamma^{-1}(mI+B)\, \Gamma(mI+B, x) \frac{v^{m}}{m!}.\label{e46}
\end{align}
\end{remark}
\begin{theorem}
The generalized incomplete exponential matrix function ${_{p}E_{q}}(x, A, B;v)$  satisfies the following integral representation:
\begin{align}
_{p}E_{q}\left[(x, A, B; v)\vert\begin{array}{c}E_1, E_2, \cdots, E_{p}\\ F_1, F_2, \cdots, F_{q}\end{array}\right]= \int_{x}^{\infty} t^{B-I} e^{-t}\,  _{p}R_{q}\left[\begin{array}{c}E_1, E_2, \cdots, E_{p}\\ F_1, F_2, \cdots, F_{q}\end{array}\vert  A, B; v t^{A}\right] dt, \label{e47}
\end{align}
where $A$, $B$, $E_i$ and $F_{j}$, $1\leq i\leq  p$, $1\leq j\leq  q$ be matrices in $\mathbb{C}^{r\times r}$ such that $F_{j}+kI$ are invertible for all integers $k\geq 0$, and $B$ is positive stable.
\end{theorem}
\begin{proof}
Using the integral representation of the incomplete gamma matrix function defined by (\ref{1eq5}), we obtain
\begin{align}
&_{p}E_{q}\left[(x, A, B; v)\vert\begin{array}{c}E_1, E_2, \cdots, E_{p}\\ F_1, F_2, \cdots, F_{q}\end{array}\right]\notag\\
&=\int_{x}^{\infty} t^{mA+B-I} e^{-t}\left( \sum_{m=0}^{\infty}\Gamma^{-1}(mA+B) (E_1)_{m}\cdots(E_p)_{m}(F_1)^{-1}_{m}\cdots(F_q)^{-1}_{m}\frac{v^{m}}{m!}\right) dt.
\end{align}
Reversing the order of summation and integration yields the R.H.S. of assertion  (\ref{e47}).
\end{proof}
\begin{corollary}
Putting   $A=I, B=C$, $p=1, q=0$ i.e. $E_{1}=A$, (\ref{e47}) reduces to 
\begin{align}
_{1}E_{0}\left[(x, I, C; v)\vert\begin{array}{c}A\\ -\end{array}\right]= \Gamma^{-1}(C)\int_{x}^{\infty} t^{C-I} e^{-t}\, _{1}F_{1}\left[\begin{array}{c}A\\ C\end{array}\vert vt\right] dt,
\end{align}
where $A$ and $C$ be matrices in $\mathbb{C}^{r\times r}$ such that $C+kI$ is invertible for all integers $k\geq 0$ and $C$ is positive stable.
\end{corollary}
\begin{corollary}
For the matrix function $_{p}R_{q}(A, B; z)$, the following integral representation holds true:
\begin{align}
&_{p}R_{q}\left[\begin{array}{c}E_1, E_2, \cdots, E_{p}\\ F_1, F_2, \cdots, F_{q}\end{array}\vert  A, B; v \right]\notag\\
&= \Gamma^{-1}(E_1)\int_{0}^{\infty} t^{E_{1}-I} e^{-t}\,  _{p-1}R_{q}\left[\begin{array}{c} E_2, \cdots, E_{p}\\ F_1, F_2, \cdots, F_{q}\end{array}\vert  A, B; v t\right]dt,
\end{align}
where $A$, $B$, $E_i$ and $F_{j}$, $1\leq i\leq  p$, $1\leq j\leq  q$ be matrices in $\mathbb{C}^{r\times r}$ such that $F_{j}+kI$ are invertible for all integers $k\geq 0$ and $E_1$ is positive stable matrix.
\end{corollary}
\begin{theorem} Let $A$, $B$, $E_i$ and $F_{j}$, $1\leq i\leq  p$, $1\leq j\leq  q$ be matrices in $\mathbb{C}^{r\times r}$ such that  $E_p F_j= F_j E_p$, $F_{j}+kI$ are invertible for all integers $k\geq 0$ and $E_1, F_1, F_1-E_1$ are positive stable.  Then 
the  matrix function ${_{p}E_{q}}(x, A, B;v)$  defined in
$(\ref{e43})$ can be put in the integral form as
\begin{align}
&_{p}E_{q}\left[(x, A, B; v)\vert\begin{array}{c}E_1, E_2, \cdots, E_{p}\\ F_1, F_2, \cdots, F_{q}\end{array}\right]\notag\\&=\int_{0}^{1} \, _{p-1}E_{q-1}\left[(x, A, B; tv)\vert\begin{array}{c} E_1, \cdots, E_{p-1}\\  F_1, \cdots, F_{q-1}\end{array}\right] t^{E_p-1}(1-t)^{F_q-E_p-1}dt\notag\\
&\times {\left[ \mathfrak{B}(E_p, F_{q}-E_{p})\right]}^{-1}.\label{4.10}
\end{align}
\end{theorem}
\begin{proof}
\begin{align}
&_{p}E_{q}\left[(x, A, B; v)\vert\begin{array}{c}E_1, E_2, \cdots, E_{p}\\ F_1, F_2, \cdots, F_{q}\end{array}\right]\notag\\&
= \sum_{m=0}^{\infty}\Gamma^{-1}(mA+B)\Gamma(mA+B, x) (E_{1})_{m}\dots (E_p)_{m} (F_{1})^{-1}_{m}\dots (F_q)^{-1}_{m}\frac{v^{m}}{m!}\notag\\
&= \sum_{m=0}^{\infty}\Gamma^{-1}(mA+B)\Gamma(mA+B, x) (E_{1})_{m}\dots (E_{p-1})_{m} (F_{1})^{-1}_{m}\dots (F_{q-1})^{-1}_{m}\frac{v^{m}}{m!}\notag\\
&\times \mathbb{B}{(E_p+mI, F_{q}-E_{p})}\, {\left[ \mathfrak{B}(E_p, F_{q}-E_{p})\right]}^{-1}\label{4.11}
\end{align}
Applying the integral representation of  beta matrix function in (\ref{4.11}), we get (\ref{4.10}). This completes the proof of this theorem. 
\end{proof}
\begin{corollary} For the matrix function $_{p}R_{q}(A, B; v)$, we have the following integral representation:
\begin{align}
&_{p}R_{q}\left[\begin{array}{c}E_1, E_2, \cdots, E_{p}\\ F_1, F_2, \cdots, F_{q}\end{array}\vert A, B; v \right]\notag\\&=\int_{0}^{1} \, _{p-1}R_{q-1}\left[\begin{array}{c} E_2, \cdots, E_{p}\\  F_2, \cdots, F_{q}\end{array}\vert A, B; vt \right] t^{E_p-1}(1-t)^{F_q-E_p-1}dt\notag\\
&\times {\left[ \mathfrak{B}(E_p, F_{q}-E_{p})\right]}^{-1},
\end{align}
where  $A$, $B$, $E_i$ and $F_{j}$, $1\leq i\leq  p$, $1\leq j\leq  q$ be matrices in $\mathbb{C}^{r\times r}$ such that  $E_p F_j= F_j E_p$, $F_{j}+kI$ are invertible for all integers $k\geq 0$ and $E_p, F_q, F_q-E_p$ are positive stable.
\end{corollary}
\begin{theorem}
Let $A$, $B$, $E_i$ and $F_{j}$, $1\leq i\leq  p$, $1\leq j\leq  q$ be matrices in $\mathbb{C}^{r\times r}$ such that  $E_i F_j= F_j E_i$, $F_{j}+kI$ are invertible for all integers $k\geq 0$ .  Then 
the generalized incomplete exponential matrix function ${_{p}E_{q}}(x, A, B;v)$  have the following derivative formula:
\begin{align}
&\frac{d^{n}}{dv^{n}}\left(_{p}E_{q}\left[(x, A, B; v)\vert\begin{array}{c}E_1, E_2, \cdots, E_{p}\\ F_1, F_2, \cdots, F_{q}\end{array}\right]\right)\notag\\&=\, _{p}E_{q}\left[(x, A, A+B; v)\vert\begin{array}{c}E_1, E_2, \cdots, E_{p}\\ F_1, F_2, \cdots, F_{q}\end{array}\right]\times (E_1)_{n}\cdots(E_p)_{n} (F_1)^{-1}_{n}\cdots(F_q)^{-1}_{n}.\label{4.13}
\end{align}
\end{theorem}
\begin{proof}
Differentiating (\ref{e43}) with repect to $v$ and replacing $m\rightarrow m+1$, we get
\begin{align}
&_{p}E_{q}\left[(x, A, B; v)\vert\begin{array}{c}E_1, E_2, \cdots, E_{p}\\ F_1, F_2, \cdots, F_{q}\end{array}\right]\notag\\&
= \sum_{m=0}^{\infty}\Gamma^{-1}(mA+A+B)\Gamma(mA+A+B, x) (E_{1})_{m+1}\dots (E_p)_{m+1} (F_{1})^{-1}_{m+1}\dots (F_q)^{-1}_{m+1}\frac{v^{m}}{m!}\notag
\end{align}
Using the relation  $(A)_{m+1}= A(A+I)_{m}$, we arrive at
\begin{align}
&\frac{d}{dv}\left(_{p}E_{q}\left[(x, A, B; v)\vert\begin{array}{c}E_1, E_2, \cdots, E_{p}\\ F_1, F_2, \cdots, F_{q}\end{array}\right]\right)\notag\\&=\, _{p}E_{q}\left[(x, A, A+B; v)\vert\begin{array}{c}E_1, E_2, \cdots, E_{p}\\ F_1, F_2, \cdots, F_{q}\end{array}\right]\times (E_1)_{1}\cdots(E_p)_{1}\, (F_1)^{-1}_{1}\cdots(F_q)^{-1}_{1}.
\end{align}
by repeating above procedure $n-$times yields the R.H.S. of assertion (\ref{4.13}).
\end{proof}
\begin{corollary}For the matrix function $_{p}R_{q}(A, B; v)$, we have the following derivative formula:
\begin{align}
&\frac{d^{n}}{dv^{n}}\left(_{p}R_{q}\left[\begin{array}{c}E_1, E_2, \cdots, E_{p}\\ F_1, F_2, \cdots, F_{q}\end{array}\vert A, B; v \right]\right)\notag\\&=\, _{p}R_{q}\left[\begin{array}{c}E_1, E_2, \cdots, E_{p}\\ F_1, F_2, \cdots, F_{q}\end{array}\vert A, A+B; v \right]\times (E_1)_{n}\cdots(E_p)_{n} (F_1)^{-1}_{n}\cdots(F_q)^{-1}_{n},
\end{align}
where  $A$, $B$, $E_i$ and $F_{j}$, $1\leq i\leq  p$, $1\leq j\leq  q$ be matrices in $\mathbb{C}^{r\times r}$ such that  $E_i F_j= F_j E_i$, $F_{j}+kI$ are invertible for all integers $k\geq 0$ .
\end{corollary}
\begin{theorem}
The generalized incomplete exponential matrix function ${_{p}E_{q}}(x, A, B;v)$  have the following  partial derivatives holds true:
\begin{align}
&\frac{\partial}{\partial v}\left(_{p}E_{q}\left[(x, A, B; v)\vert\begin{array}{c}E_1, E_2, \cdots, E_{p}\\ F_1, F_2, \cdots, F_{q}\end{array}\right]\right)\notag\\&= \, _{p}E_{q}\left[(x, A, A+ B; v)\vert\begin{array}{c}E_1+I, E_2+I, \cdots, E_{p}+I\\ F_1+I, F_2+I, \cdots, F_{q}+I\end{array}\right]\notag\\&\times  E_1\cdots E_{p}\,  F_{1}^{-1}\cdots F_{q}^{-1}, \label{e4.16}
\end{align}
\begin{align}
&\frac{\partial}{\partial x}\left(_{p}E_{q}\left[(x, A, B; v)\vert\begin{array}{c}E_1, E_2, \cdots, E_{p}\\ F_1, F_2, \cdots, F_{q}\end{array}\right]\right)\notag\\&= - e^{-x} x^{B-I} \left(_{p}R_{q}\left[\begin{array}{c}E_1, E_2, \cdots, E_{p}\\ F_1, F_2, \cdots, F_{q}\end{array}\vert A, B; vx^{A} \right]\right), \label{4.17}
\end{align}
where  $A$, $B$, $E_i$ and $F_{j}$, $1\leq i\leq  p$, $1\leq j\leq  q$ be matrices in $\mathbb{C}^{r\times r}$ such that  $E_i F_j= F_j E_i$, $F_{j}+kI$ are invertible for all integers $k\geq 0$.
\end{theorem}
\begin{proof}
Differentiating partially (\ref{e43}) with respect to $v$, we get
\begin{align}
&\frac{\partial}{\partial v}\left(\,_{p}E_{q}\left[(x, A, B; v)\vert\begin{array}{c}E_1, E_2, \cdots, E_{p}\\ F_1, F_2, \cdots, F_{q}\end{array}\right]\right)\notag\\&= \,\frac{\partial}{\partial v}\left(\sum_{m=0}^{\infty}\Gamma^{-1}(mA+B)\Gamma(mA+B, x) (E_{1})_{m}\dots (E_p)_{m} (F_{1})^{-1}_{m}\dots (F_q)^{-1}_{m}\frac{v^{m}}{m!}\right), \notag\\
&=\sum_{m=1}^{\infty}\Gamma^{-1}(mA+B)\Gamma(mA+B, x) (E_{1})_{m}\dots (E_p)_{m} (F_{1})^{-1}_{m}\dots (F_q)^{-1}_{m}\frac{v^{m-1}}{(m-1)!}.\notag
\end{align}
This leads to proof of (\ref{e4.16}) by replacing $m$$\rightarrow$$ m+1$.\\

We differentiate partially (\ref{e47}) with respect to $x$ to demonstrate (\ref{4.17}).
\end{proof}
\begin{theorem}
For the generalized incomplete exponential matrix function ${_{p}E_{q}}(x, A, B;v)$, the following addition formula holds true:  
\begin{align}
& _{p}E_{q}\left[(x, A, B; w+v)\vert\begin{array}{c}E_1, E_2, \cdots, E_{p}\\ F_1, F_2, \cdots, F_{q}\end{array}\right]\notag\\
&= \sum_{n=0}^{\infty}\, _{p}E_{q}\left[(x, A, nA+B; w)\vert\begin{array}{c}E_1+nI, E_2+nI, \cdots, E_{p}+nI\\ F_1+nI, F_2+nI, \cdots, F_{q}+nI\end{array}\right] \frac{v^{n}}{n!}\notag\\
&\times \Gamma(E_1+nI)\cdots \Gamma(E_{p}+nI)\Gamma(F_1+nI)^{-1}\cdots \Gamma(F_{q}+nI)^{-1}\notag\\
&\times  \Gamma(E_1)^{-1}\cdots \Gamma(E_{p})^{-1}\Gamma(F_1)\cdots \Gamma(F_{q}),\label{4.18}
\end{align}
where  $A$, $B$, $E_i$ and $F_{j}$, $1\leq i\leq  p$, $1\leq j\leq  q$ be matrices in $\mathbb{C}^{r\times r}$ such that  $E_i F_j= F_j E_i$, $F_{j}+kI$ are invertible for all integers $k\geq 0$.
\end{theorem}
\begin{proof}
Applying the definition of generalized incomplete exponential matrix function ${_{p}E_{q}}(x, A, B; v)$ given in (\ref{e43}), expanding L.H.S. of (\ref{4.18}) and using the identity \cite{hm1}
\begin{align}
\sum_{N=0}^{\infty} f(N) \frac{(w+v)^{N}}{N!}= \sum_{m, n=0}^{\infty} f(m+n) \frac{w^{m}}{m!}\frac{v^{n}}{n!}.\notag
\end{align}
We get (\ref{4.18}). This finishes the proof of this theorem.
\end{proof}
\begin{theorem}
For the generalized incomplete exponential matrix function ${_{p}E_{q}}(x, A, B;v)$, the following multiplication formula holds true:  
\begin{align}
& _{p}E_{q}\left[(x, A, B; wv)\vert\begin{array}{c}E_1, E_2, \cdots, E_{p}\\ F_1, F_2, \cdots, F_{q}\end{array}\right]\notag\\
&= \sum_{n=0}^{\infty}\, _{p}E_{q}\left[(x, A, nA+B; w)\vert\begin{array}{c}E_1+nI, E_2+nI, \cdots, E_{p}+nI\\ F_1+nI, F_2+nI, \cdots, F_{q}+nI\end{array}\right] \frac{w^{n}(v-1)^{n}}{n!}\notag\\
&\times \Gamma(E_1+nI)\cdots \Gamma(E_{p}+nI)\Gamma(F_1+nI)^{-1}\cdots \Gamma(F_{q}+nI)^{-1}\notag\\
&\times  \Gamma(E_1)^{-1}\cdots \Gamma(E_{p})^{-1}\Gamma(F_1)\cdots \Gamma(F_{q}),\label{4.19}
\end{align}
where  $A$, $B$, $E_i$ and $F_{j}$, $1\leq i\leq  p$, $1\leq j\leq  q$ be matrices in $\mathbb{C}^{r\times r}$ such that  $E_i F_j= F_j E_i$, $F_{j}+kI$ are invertible for all integers $k\geq 0$.
\end{theorem}
\begin{proof}
Proof of above theorem is similar to  theorem $4.5$.
\end{proof}
\begin{theorem}
For the generalized incomplete exponential matrix function ${_{p}E_{q}}(x, A, B;v)$, we have the following integral representation:  
\begin{align}
&\int_{0}^{t} v^{A-I}(t-v)^{B-I} \, _{p}E_{q}\left[(x, A, B; \lambda v^{k})\vert\begin{array}{c}E_1, E_2, \cdots, E_{p}\\ F_1, F_2, \cdots, F_{q}\end{array}\right] dv\notag\\
&= \mathfrak{B}(A, B) \,t^{A+B-I} \, _{p+k}E_{q+k}\left[(x, A, B; \lambda t^{k})\vert\begin{array}{c}\Delta (k, A), E_1, E_2, \cdots, E_{p}\\\Delta (k, A+B),  F_1, F_2, \cdots, F_{q}\end{array}\right],\label{4.20} 
\end{align}
where  $A$, $B$, $E_i$ and $F_{j}$, $1\leq i\leq  p$, $1\leq j\leq  q$ be matrices in $\mathbb{C}^{r\times r}$ such that   $F_{j}+kI$ are invertible for all integers $k\geq 0$. Also $ \Delta (k, A)$ stands for $\frac{A}{k}$, $ \frac{A+I}{k}$, \dots, $\frac{A+(k-1)I}{k}$.
\end{theorem}
\begin{proof}
Suppose $\mathfrak{L}$ be the L.H.S. of (\ref{4.20}). Then, using (\ref{e43}), this gives
\begin{align}
\mathfrak{L}&= \int_{0}^{t} v^{A-I} (t-v)^{B-I}\sum_{m=0}^{\infty}\Gamma^{-1}(mA+B) \Gamma(mA+B, x)\notag\\&\times (E_1)_{m}\cdots (E_p)_{m} (F_1)^{-1}_{m}\cdots (F_q)^{-1}_{m} \frac{(\lambda v^{k})^{m}}{m!} dv.\notag
\end{align}
Putting $v=t y$, we get
\begin{align}
\mathfrak{L}&=  t^{A+B-I}\int_{0}^{1}\sum_{m=0}^{\infty} y^{A+(km-1)I} (1-y)^{B-I}\,\Gamma^{-1}(mA+B) \Gamma(mA+B, x)\notag\\&\times (E_1)_{m}\cdots (E_p)_{m} (F_1)^{-1}_{m}\cdots (F_q)^{-1}_{m} \frac{(\lambda t^{k})^{m}}{m!} dv,\notag\\
&= t^{A+B-I}\sum_{m=0}^{\infty}\mathfrak{B}(A+kmI, B) \Gamma^{-1}(mA+B) \Gamma(mA+B, x)\notag\\&\times (E_1)_{m}\cdots (E_p)_{m} (F_1)^{-1}_{m}\cdots (F_q)^{-1}_{m} \frac{(\lambda t^{k})^{m}}{m!} dv,\notag\\
&=  t^{A+B-I}\sum_{m=0}^{\infty}\Gamma(A+kmI)\Gamma(B)\Gamma^{-1}(A+B+kmI) \Gamma^{-1}(mA+B) \Gamma(mA+B, x)\notag\\&\times (E_1)_{m}\cdots (E_p)_{m} (F_1)^{-1}_{m}\cdots (F_q)^{-1}_{m} \frac{(\lambda t^{k})^{m}}{m!} dv.\notag
\end{align}
Now using the property of Pochhammer  matrix symbol defined in (\ref{g1}), this leads to the R.H.S. of (\ref{4.20}).
\end{proof}
\begin{theorem} Let $A$, $B$, $C$, $E_i$ and $F_{j}$, $1\leq i\leq  p$, $1\leq j\leq  q$ be matrices in $\mathbb{C}^{r\times r}$ such that   $F_{j}+kI$ are invertible for all integers $k\geq 0$ and $B$, $C$, $C+B$ are positive stable. Then  the generalized incomplete exponential matrix function ${_{p}E_{q}}(x, A, B;v)$, satisfies the following integral representation:  
\begin{align}
&\int_{t}^{y} (y-v)^{C-I}(v-t)^{B-I} \, _{p}E_{q}\left[(x, A, B; \lambda (v-t)^{k})\vert\begin{array}{c}E_1, E_2, \cdots, E_{p}\\ F_1, F_2, \cdots, F_{q}\end{array}\right] dv\notag\\
&= \mathfrak{B}(B, C) \,(y-t)^{C+B-I} \, _{p+k}E_{q+k}\left[(x, A, B; \lambda (x-t)^{k})\vert\begin{array}{c}\Delta (k, B), E_1, E_2, \cdots, E_{p}\\\Delta (k, B+C),  F_1, F_2, \cdots, F_{q}\end{array}\right],\label{4.21} 
\end{align}
where $ \Delta (k, A)$ stands for $\frac{A}{k}$, $ \frac{A+I}{k}$, \dots, $\frac{A+(k-1)I}{k}$.
\end{theorem}
\begin{proof}
Suppose $\mathfrak{A}$ be the L.H.S. of (\ref{4.21}). Then, using (\ref{e43}), this gives
\begin{align}
\mathfrak{A}&= \int_{t}^{y} (y-v)^{C-I}(v-t)^{B-I}\sum_{m=0}^{\infty}\Gamma^{-1}(mA+B) \Gamma(mA+B, x)\notag\\&\times (E_1)_{m}\cdots (E_p)_{m} (F_1)^{-1}_{m}\cdots (F_q)^{-1}_{m} \frac{(\lambda (v-t)^{k})^{m}}{m!} dv.\notag
\end{align}
Putting $s=\frac{v-t}{y-t}$, we get
\begin{align}
\mathfrak{A}&=  (y-t)^{C+B-I}\int_{0}^{1}\sum_{m=0}^{\infty} s^{B+(km-1)I} (1-s)^{C-I}\,\Gamma^{-1}(mA+B) \Gamma(mA+B, x)\notag\\&\times (E_1)_{m}\cdots (E_p)_{m} (F_1)^{-1}_{m}\cdots (F_q)^{-1}_{m} \frac{(\lambda (y-t)^{k})^{m}}{m!} dv,\notag\\
&= (y-t)^{C+B-I}\sum_{m=0}^{\infty}\mathfrak{B}(B+kmI, C) \Gamma^{-1}(mA+B) \Gamma(mA+B, x)\notag\\&\times (E_1)_{m}\cdots (E_p)_{m} (F_1)^{-1}_{m}\cdots (F_q)^{-1}_{m} \frac{(\lambda (y-t)^{k})^{m}}{m!} dv,\notag\\
&=  (y-t)^{C+B-I}\sum_{m=0}^{\infty}\Gamma(B+kmI)\Gamma(C)\Gamma^{-1}(B+C+kmI) \Gamma^{-1}(mA+B) \Gamma(mA+B, x)\notag\\&\times (E_1)_{m}\cdots (E_p)_{m} (F_1)^{-1}_{m}\cdots (F_q)^{-1}_{m} \frac{(\lambda (y-t)^{k})^{m}}{m!} dv.\notag
\end{align}
Now using the property of Pochhammer  matrix symbol defined in (\ref{g1}), this leads to the R.H.S. of (\ref{4.21}).
\end{proof}
\begin{theorem}
The  incomplete exponential matrix function ${_{2}E_{1}}(x, A, B; v)$ satisfy the following relaton
\begin{align}
&_{2}E_{1}\left[(x, A, B; 1)\vert\begin{array}{c}E_1, E_2\\ F_1\end{array}\right]\notag\\&= \Gamma(F_{1}-E_{1}-E_{2})\Gamma(F_1)\Gamma^{-1}(F_{1}-E_{1})\Gamma^{-1}(F_{1}-E_{2})-\gamma(mA+B, x)\times {_{2}}R_{1}(A, B;1),\label{4.22}
\end{align}
where  $A$, $B$, $E_1$, $E_2$ and $F_{1}$  be matrices in $\mathbb{C}^{r\times r}$ such that $F_{1}+kI$  is  invertible for all integers $k\geq 0$,  and $F_1$, $F_{1}-E_{1}$, $F_{1}-E_{2}$  and $F_{1}-E_{1}-E_{2}$ are positive stable matrices where all matrices are commutative.
\end{theorem}
\begin{proof}
Putting $p=2$, $q=1$, $v=1$ in  (\ref{e44}), we obtain
\begin{align}
&_{2}E_{1}\left[(x, A, B; 1)\vert\begin{array}{c}E_1, E_2\\ F_1\end{array}\right]\notag\\&= {_{2}}F_{1}\left[\begin{array}{c}E_1, E_2\\ F_1\end{array}\vert 1\right]- {_{2}}e_{1}\left[(x, A, B; 1)\vert\begin{array}{c}E_1, E_2\\ F_1\end{array}\right]\notag\\
&= {_{2}}F_{1}\left[\begin{array}{c}E_1, E_2\\ F_1\end{array}\vert 1\right]- \int_{0}^{x} t^{B-I} e^{-t}{_{2}}R_{1}\left[\begin{array}{c}E_1, E_2\\ F_1\end{array}\vert A, B; t^{A}\right]dt,\notag\\
&= {_{2}}F_{1}\left[\begin{array}{c}E_1, E_2\\ F_1\end{array}\vert 1\right]- \int_{0}^{x} t^{B-I} e^{-t}\sum_{m=0}^{\infty}\Gamma^{-1}(mA+B) (E_1)_{m} (E_{2})_{m} (F_{1})^{-1}_{m}\frac{(t^{A})^{m}}{m!} dt.
\end{align}
Applying  Gauss summation matrix formula \cite{As} for v = 1 and changing the order of summation and integration, we get
\begin{align}
&_{2}E_{1}\left[(x, A, B; 1)\vert\begin{array}{c}E_1, E_2\\ F_1\end{array}\right]\notag\\&=  \Gamma(F_{1}-E_{1}-E_{2})\Gamma(F_1)\Gamma^{-1}(F_{1}-E_{1})\Gamma^{-1}(F_{1}-E_{2})\notag\\&- \sum_{m=0}^{\infty} \Gamma^{-1}(mA+B) (E_1)_{m} (E_{2})_{m} (F_{1})^{-1}_{m}\frac{1}{m!}\times\left( \int_{0}^{x} t^{mA+B-I} e^{-t}dt\right).
\end{align}
After some  simplification by using (\ref{1eq4})and (\ref{9}), we obtain R.H.S. of (\ref{4.22}).
\end{proof}
\begin{theorem} Let $A$, $B$, $E_1$, $E_2$ and $F_{1}$  be matrices in $\mathbb{C}^{r\times r}$ such that  $E_1 F_1= F_1 E_1$, $E_{1}E_{2}= E_{2}E_{1}$,  $F_{1}+kI$  is  invertible for all integers $k\geq 0$ and $E_1$, $E_1-F_1+I$, $ F_1-I$ are positive stable. Then the  matrix function ${_{2}E_{1}}(x, A, B;v)$  satisfies the  following recurrence relation:
\begin{align}
&_{2}E_{1}\left[(x, A, B; v)\vert\begin{array}{c}E_1, E_2\\ F_1\end{array}\right] (E_{1}-F_{1}+I)\notag\\
&= \, _{2}E_{1}\left[(x, A, B; v)\vert\begin{array}{c}E_1+I, E_2\\ F_1\end{array}\right] E_{1}-\,_{2}E_{1}\left[(x, A, B; v)\vert\begin{array}{c}E_1, E_2\\ F_1-I\end{array}\right] (F_{1}-I).\label{4.16}
\end{align}
\end{theorem}
\begin{proof}
Using the definition of the incomplete exponential matrix function $_2E _1(x, A, B;v)$, take $p=2, q=1$ in (\ref{e43}), and the R.H.S. of  (\ref{4.16}) is of the form
\begin{align}
R.H.S.&= \sum_{m=0}^{\infty}\Gamma^{-1}(mA+B)\Gamma(mA+B, x) (E_1+I)_{m} (E_2)_{m}(F_{1})^{-1}_{m} \frac{v^m}{m!}\, E_{1}\notag\\
&-\sum_{m=0}^{\infty}\Gamma^{-1}(mA+B)\Gamma(mA+B, x) (E_1)_{m} (E_2)_{m}(F_{1}-I)^{-1}_{m} \frac{v^m}{m!}\, (F_{1}-I)
\end{align}
 Using  the   relations
\begin{align}
E_1 (E_1+I)_{m}&= (E_{1}+mI) (E_1)_{m},\notag\\
(F_1-I) (F_{1}-I)_{m}^{-1}&= (F_{1}+(m-1)I) (F_{1})_{m}^{-1},\notag
\end{align}
this yields the L.H.S. of (\ref{4.16}).
\end{proof}
\section{Conclusion}
 We have investigated the  generalized  incomplete exponential matrix functions  ${_{p}e_{q}}(x, A, B;v)$ and ${_{p}E_{q}}(x, A, B;v)$.  We have found several characteristics of each of these  generalized  incomplete exponential matrix functions, for example  various integral representations, addition  formulas and derivative formulas etc.  Some integral representations containing incomplete gamma matrix function,  Bessel and modified Bessel matrix  functions are also presented. Inference of this article is generalization of  \cite{ APa, AQ} for the matrix cases.  These generalized incomplete exponential matrix functions have applications in many areas of mathematics and mathematical physics.\\

{\bf Acknowledgments.}  The second author is grateful to the University Grants Commission of India for financial assistance in the form of a Junior Research Fellowship.\\

{\bf Declarations}\\

{\bf Conflict of interest} The authors declare no conflict of interest.  

\end{document}